\def\BE{\begin{equation}}
\def\EE#1{\label{#1}\end{equation}}

\def\T{{\mathbb T}}
\def\subs#1#2{\mbox{\small${#1\atop#2}$}}
\def\nL{\nabla^{\rm L}}
\def\rf#1{(\ref{#1})}
\def\bq{{\bm q}}
\def\bv{{\bm v}}
\def\bx{{\bm x}}
\def\bxi{{\bm\xi}}
\def\ba{{\bm a}}
\def\bom{{\bm\omega}}
\documentclass[cmp]{svjour}
\usepackage{bm,amsmath,amssymb,hyperref}
\begin{document}
\title{A very smooth ride in a rough sea}
\author{Uriel Frisch\inst{1}\and Vladislav Zheligovsky\inst{2,1}}
\institute{UNS, CNRS, Lab. Lagrange, OCA, B.P. 4229, 06304 Nice Cedex 4, France
\and Institute of earthquake prediction theory and mathematical geophysics
of the Russ.~Ac.~Sci., 84/32 Profsoyuznaya St, 117997 Moscow, Russian Federation}

\maketitle
\begin{abstract}
It has been known for some time that a 3D incompressible Euler flow
that has initially a barely smooth velocity field nonetheless has
Lagrangian fluid particle trajectories that are analytic in time for
at least a finite time (Ph.~Serfati C.~R.~Acad.~Sci. S\'erie I
{\bf 320}, 175--180 (1995); A.~Shnirelman arXiv:1205.5837 (2012)).
Here an elementary derivation is given, based on Cauchy's form of the Euler
equations in Lagrangian coordinates. This form implies simple
recurrence relations among the time-Taylor coefficients of the
Lagrangian map, used here to derive bounds for the C$^{1,\gamma}$
H\"older norms of the coefficients and infer temporal analyticity of
Lagrangian trajectories when the initial velocity is C$^{1,\gamma}$.
\end{abstract}

\section{Introduction}

The issue of finite time blow-up of solutions to the 3D Euler equation
for incompressible fluid is still an open question. The available results
indicate that (a hypothetical) occurrence of singularity is intimately related
to the loss of spatial regularity (in the Eulerian coordinates) of the solution.
It is known, however, that an initial regularity of the flow, that is quite modest ---
marginally better than $C^1$, which classical solutions to the Euler equations
ought to possess --- will preserve this regularity for at least a finite time
$t_c$ (see, e.g., \cite{BaFr}). It may~be surprising to find that,
in such a rough velocity field, the time dependence of the position of any
fluid particle, i.e. the Lagrangian time dependence, is actually analytic up to
$t_c$. By contrast, in Eulerian coordinates, lack of spatial smoothness also
translates into lack of temporal smoothness, because of the sweeping of
fine-scale structure by large-scale velocity fields. For instance, the turbulent
solutions \`a la Kolmogorov--Onsager have a spatial regularity of the
H\"older type with an exponent close to $1/3$ in both the spatial and the
temporal Eulerian domain.

The unexpected possibility that Lagrangian trajectories could be
analytic in time was first considered by Serfati \cite{Se92} (part
2), \cite{Se95a,Se95b} (see also \cite{Ga94}) who suggested to
attack this problem as an abstract ordinary differential equation
(ODE) in a complex Banach space with an analytic right hand
side. Prior to this, Chemin \cite{Ch92} proved that the Lagrangian
trajectories are $C^\infty$ smooth. More recently, Shnirelman \cite{Shni},
using the geometric interpretation of the Euler equation, obtained another
analytic abstract ODE formulation; he then constructed a time-analytic
Lagrangian solution to the ODE applying Picard's theorem.
The fact, that the Lagrangian structure can be nicer than the Eulerian
one, was known to Ebin and Marsden \cite{EM70} who observed
(on p.~104) that in Lagrangian coordinates the Euler equations can be
written in such a way that ``no derivative loss occurs''. An early
embodyment of this is actually found already in Cauchy's little-known
Lagrangian formulation of the Euler equation \cite{Cau27}.

This formulation will be our starting point. It has the advantage
of allowing to set up not only an elementary proof of the analyticity, but also
a procedure for calculation of the various temporal Taylor coefficients.
In Section~\ref{s:basic}, we derive simple recursion relations for the
time-Taylor coefficients of the Lagrangian displacement. Then, in
Section~\ref{s:hoelder}, we obtain bounds for the C$^{1,\gamma}$ of
these coefficients. The key to obtaining these bound is that the
Lagrangian gradient of $n$-th Taylor coefficient can be recursively
written in terms of nonlocal operators applied to products of two and
three gradients of lower-order such Taylor coefficients. These
nonlocal operators are essentially inverse Laplacians composed with
two suitable space derivatives and thus all the Taylor coefficients
stay within the same H\"older space and simple estimates are obtained
for their H\"older norms. Concluding remarks are presented in
Section~\ref{s:conclusion}.

\section{Basic equations ...}\label{s:basic}

Our starting point is not the usual Eulerian formulation of the
Euler equations for incompressible fluid, but a little-known Lagrangian
formulation due to Cauchy \cite{Cau27}, specifically his eq.~(15), together
with the condition of unit Jacobian of the Lagrangian map that expresses
incompressibility. We use $\bq$ and $\bx$ for the Lagrangian and Eulerian
positions of fluid particles and $\bxi\equiv\bx-\bq$ for the displacement.
In modern notation the dynamical equations of Cauchy \cite{Cau27} take the form
\begin{eqnarray*}
\sum_{k=1}^3\nL\dot{x}_k\times\nL x_k&=&\bom_0,\\
\det\,(\nL\bx)&=&1,\\
\end{eqnarray*}
where $\bom_0\equiv\nL\times\bv_0$ is the initial vorticity, $\nL$ denotes
the gradient in the Lagrangian variables and $\nL\bx$ the Jacobian matrix
with entries $\nL_i x_j$. The Cauchy equations can be obtained as the (Lagrangian) curl
of the Weber equations \cite{We68} (see also \cite{Const00,bennett,ohkitani}). These equations,
obviously, express the time invariance of the l.h.s. and are therefore sometimes
referred to as ``Cauchy invariants of an ideal fluid''.
For simplicity, in what follows we assume spatial $2\pi$-periodicity (flow
in the 3D torus $\T^3=[0,2\pi]^3$), but this can be relaxed.

Following \cite{FM12}, we study the formal Taylor expansion
of the particle displacement $\bxi(\bq,t)$ in a power series in time,
\BE\bxi(\bq,t)=\sum_{s=1}^\infty\bxi^{(s)}(\bq)t^s.\EE{xiser}
We shall introduce a power
series in $t$ with real non-negative coefficients bounding certain norms
of $\bxi^{(s)}(\bq)$, and prove, under the appropriate assumptions, that this series
converges for a certain real $t=t_c$. We shall conclude that, for complex~$t$,
series \rf{xiser} has a finite radius of convergence, not smaller than $t_c$.

In terms of the displacement, Cauchy's equations become:
\begin{eqnarray}
&&\nL\times\dot{\bxi}+\sum_{k=1}^3\nL\dot{\xi}_k\times\nL\xi_k=\bom_0,
\label{3DDisplacement1}\\
&&\nL\cdot\bxi+\sum_{1\le i<j\le3}\left(\rule{0mm}{1em}
(\nL_i\xi_i)\nL_j\xi_j-(\nL_i\xi_j)\nL_j\xi_i\right)+\det(\nL\bxi)=0.
\label{3DDisplacement2}\end{eqnarray}
Following \cite{FM12}, we now derive recurrence relations for
the time-Taylor coefficients of the displacement at $t=0$. Substituting
\rf{xiser} into \rf{3DDisplacement1} and \rf{3DDisplacement2}, we find
\begin{eqnarray}
s\,\nL\times\bxi^{(s)}&=&\bom_0\delta^s_1-\!\!
\sum_{\subs{k=1,2,3}{0<n<s}}\!\!n\nL\xi^{(n)}_k\times\nL\xi^{(s-n)}_k,\label{Frev}\\
\nL\!\cdot\bxi^{(s)}&=&\!\!\sum_{\subs{1\le i<j\le3}{0<n<s}}\!\!
\left(\rule{0mm}{1em}(\nL_j\xi^{(n)}_i)\nL_i\xi^{(s-n)}_j-(\nL_i\xi^{(n)}_i)
\nL_j\xi^{(s-n)}_j\right)\nonumber\\
&&-\!\!\sum_{\subs{i,j,k}{l+m+n=s}}\!\!
\varepsilon_{ijk}(\nL_i\xi^{(l)}_1)(\nL_j\xi^{(m)}_2)\nL_k\xi^{(n)}_3,
\label{Frst}\end{eqnarray}
where $\varepsilon_{ijk}$ is the unit antisymmetric tensor. It is immediately
seen from the two equations for $s=1$ that $\bxi^{(1)}=\bv_0$.

To handle the equations \rf{Frev}--\rf{Frst} for the Lagrangian curl and
divergence of the Taylor coefficients, we apply the Hodge decomposition,
$$\bxi^{(s)}(\bq)=\nL\times\ba^{(s)}+\nL b^{(s)},$$
where the potentials satisfy the gauge conditions: $\nL\cdot\ba^{(s)}=0$ and
vanishing of averages over the periodicity box. This yields a pair of Poisson
equations for the vector and scalar potentials $\ba^{(s)}$ and $b^{(s)}$:
\begin{eqnarray}
s\nabla^2\ba^{(s)}\!\!&=&-\bom_0\delta^s_1+\!\sum_{\subs{k=1,2,3}{0<n<s}}
n\nL\left((\nL\times\ba^{(n)})_k+\nL_kb^{(n)}\rule{0mm}{1em}\right)\nonumber\\
&&\hspace*{2em}\times\,\nL\left((\nL\!\times\ba^{(s-n)})_k+\nL_kb^{(s-n)}\right),\label{Hoev}\\
\rule{0mm}{2em}\nabla^2b^{(s)}\!\!&=&\!\!\!\!\sum_{\subs{1\le i<j\le3}{0<n<s}}\!\!
\left(\rule{0mm}{3ex}\!\nL_j\!\left(\!(\nL\!\times\ba^{(n)})_i\!+\nL_ib^{(n)}\right)
\!\nL_i\!\left(\!(\nL\!\times\ba^{(s-n)})_j\!+\nL_jb^{(s-n)}\right)\right.\nonumber\\
&&\hspace*{2em}\left.-\,\nL_i\!\left(\!(\nL\!\times\ba^{(n)})_i\!+\nL_ib^{(n)}\rule{0mm}{1em}\right)\!
\nL_j\!\left(\!(\nL\!\times\ba^{(s-n)})_j\!+\nL_jb^{(s-n)}\right)\rule{0mm}{3ex}\!\!\!\right)\nonumber\\
&&-\!\!\!\sum_{\subs{j_1,j_2,j_3}{n_1+n_2+n_3=s}}\!\!\!\varepsilon_{j_1j_2j_3}
\prod_{i=1}^3\nL_{j_i}\!\left((\nL\times\ba^{(n_i)})_i+\nL_ib^{(n_i)}\right)
\label{Host}\end{eqnarray}
(here $\nabla^2\equiv(\nL)^2$ denotes the Laplacian in the Lagrangian variables
and\break$(\nL\!\times\ba^{(n)})_k$ the $k$-th component of the vector
$\nL\!\times\ba^{(n)}$).

\section{... and their consequence}\label{s:hoelder}

In this section we establish the finite-time analyticity in time of
the Lagrangian map for solutions to the three-dimensional Euler equation.
The assumption is that the initial velocity field $\bv_0$ is in C$^{1,\gamma}$.
In other words, the initial vorticity $\bom_0\equiv\nL\times\bv_0$ is
in C$^{0,\gamma}$ for some $0<\gamma<1$, i.e., satisfies the H\"older
condition $|\bom_0(\bq+\delta\bq)-\bom_0(\bq)|\le C|\delta\bq|^\gamma$.
Demonstration of temporal analyticity amounts to proving the following statement:

\medskip
{\it Theorem.} Consider a space-periodic three-dimensional flow
of incompressible fluid governed by the Euler equation. Suppose the initial
vorticity $\bom_0(\bq)$ is in C$^{0,\gamma}(\T^3)$. Then, at sufficiently small
times, the displacement of fluid particles, $\xi(\bq,t)$, is given by power
series \rf{xiser}, whose coefficients can be recursively calculated
by solving equations \rf{Frev} and~\rf{Frst}. The radius of convergence
of \rf{xiser}, $\tau$, is bounded from below by a strictly positive quantity,
which is inversely proportional to $|\bom_0|_{0,\gamma}$. (More precisely,
$\tau\ge Q_c|\bom_0|_{0,\gamma}^{-1}/2$, where $Q_c$ is given by \rf{qc} below.)

\medskip
{\it Proof.}
From recurrence relations \rf{Hoev} and \rf{Host} we derive bounds for suitable
H\"older norms of the potentials. Using our assumption that $\bom_0$ is in
C$^{0,\gamma}$, it is easy to establish by induction that each potential
$\ba^{(s)}$ and $b^{(s)}$ is C$^{2,\gamma}$-regular. For this, we recall
that C$^{m,\gamma}(\T^3)$ is an algebra (i.e. the sum and the product
of two functions from C$^{m,\gamma}(\T^3)$ belong to this space; in particular,
$|fg|_{0,\gamma}\le|f|_{0,\gamma}|g|_{0,\gamma}$) and apply standard
inequalities for H\"older norms of functions whose integrals over
the periodicity box vanish:
\begin{eqnarray}
|\phi|_{2,\gamma}\le\widetilde{\Theta}|\nabla^2\phi|_{0,\gamma}&\qquad&
\forall\phi\in{\rm C}^{2,\gamma}(\T^3),\label{lapb1}\\
|\nL_i\nL_j\phi(\bq)|_{0,\gamma}\le\Theta|\nabla^2\phi|_{0,\gamma}&\qquad&
\forall i,j\mbox{~and~}\forall\phi\in{\rm C}^{2,\gamma}(\T^3)
\label{lapb2}\end{eqnarray}
(see, e.g., \cite{Gu27,Gu34,LaUr} or \cite{GT98}).
Here $\widetilde{\Theta}$ and $\Theta$ are some positive constants depending on
$\gamma>0$ and bounded by $C/\gamma$, $\nL_i\nL_j$ denotes the second
partial derivative $\partial^2/\partial q_i\partial q_j$, and
$|\cdot|_{s,\gamma}$ designates the H\"older norm in C$^{s,\gamma}(\T^3)$
for a scalar function $\phi$ and $|{\bm\phi}(\bq)|_{0,\gamma}=\max_k|\phi_k(\bq)|_{0,\gamma}$
for a vector-valued function $\bm\phi$.

Now, we consider the generating functions
$$A(t)\equiv\sum_{s=1}^\infty|\nabla^2\ba^{(s)}|_{0,\gamma}\,t^s,\qquad
B(t)\equiv\sum_{s=1}^\infty|\nabla^2b^{(s)}|_{0,\gamma}\,t^s,$$ thus
reintroducing time into the problem. Note that, by construction, for
$t>0$, $A(t)$ and $B(t)$ are monotonically increasing functions of time. When no
confusion is possible, we shall write just $A$ and $B$ for $A(t)$ and
$B(t)$. In order to show that the series $A$ and $B$ converge, we
proceed as follows. We derive from relations
\rf{Hoev} and \rf{Host} two inequalities for the quantities $A$ and $B$.
The first of these inequalities involves time derivatives $\dot{A}$ and
$\dot{B}$, but it can be integrated in time. The inequality resulting from this
integration and the second inequality are of the form of upper bounds for $A$
and $B$. These bounds are polynomials of $2A+B$ that involve, apart from
the contribution from the initial vorticity $\bom_0(\bq)$, only quadratic and
cubic terms. Since $\ba^{(s)}$ and $b^{(s)}$ vanish at $t=0$, the inequalities
yield an upper bound for $2A+B$ at least at small times, when $2A+B={\rm O}(t)$
and higher-order terms are insignificant.

Relation \rf{Hoev} and inequalities \rf{lapb2} imply
\begin{eqnarray*}
\dot{A}&=&\sum_{s=1}^\infty\left|-\bom_0\delta^s_1\rule{0mm}{1em}\right.
\!+\!\!\!\sum_{\subs{k=1,2,3}{0<n<s}}\!\!n
\nL\left((\nL\times\ba^{(n)})_k+\nL_kb^{(n)}\rule{0mm}{1em}\right)\\
&&\hspace*{8em}\times\left.\nL\left((\nL\!\times\ba^{(s-n)})_k+\nL_kb^{(s-n)}\rule{0mm}{1em}\right)
\right|_{0,\gamma}\!t^{s-1}\\
&\le&|\bom_0|_{0,\gamma}\!+6\!\!\sum_{\subs{s\ge1,}{0<n<s}}\!\!n
\left(\!2\max_{i,j}|\nL_i\nL_j\ba^{(n)}|_{0,\gamma}\!+\max_{i,j}|\nL_i\nL_jb^{(n)}|_{0,\gamma}\!\right)\\
&&\hspace*{6.7em}\mbox{\tiny$\times$}\left(\!2\max_{i,j}|\nL_i\nL_j\ba^{(s-n)}|_{0,\gamma}\!
+\max_{i,j}|\nL_i\nL_jb^{(s-n)}|_{0,\gamma}\!\right)t^{s-1}\\
&\le&|\bom_0|_{0,\gamma}\!+6\Theta^2\!\!\sum_{\subs{s\ge1,}{0<n<s}}\!\!n
\left(2|\nabla^2\ba^{(n)}|_{0,\gamma}\!+|\nabla^2b^{(n)}|_{0,\gamma}\!\right)\\
&&\hspace*{8em}\mbox{\tiny$\times$}\left(2|\nabla^2\ba^{(s-n)}|_{0,\gamma}\!+|\nabla^2b^{(s-n)}|_{0,\gamma}\!\right)t^{s-1}\\
&=&|\bom_0|_{0,\gamma}+6\Theta^2(2\dot{A}+\dot{B})(2A+B).
\end{eqnarray*}
Integration of this inequality in time yields
$$A\le|\bom_0|_{0,\gamma}\,t+3\Theta^2(2A+B)^2.$$
Relations \rf{Host} and inequalities \rf{lapb2} imply
$$B\le6\Theta^2(2A+B)^2+6\Theta^3(2A+B)^3.$$
From the last two inequalities we deduce
\BE p(\zeta)\equiv6\Theta^3\zeta^3+12\Theta^2\zeta^2-\zeta+Q\ge0,\EE{poly}
where we have denoted $\zeta\equiv2A+B$ and $Q\equiv2|\bom_0|_{0,\gamma}\,t$.
The discriminant of the polynomial $p(\zeta)$,
$$\Delta\equiv972\,\Theta^6\left(-Q^2-Q\left({64\over9}+{4\over3\Theta}\right)
+{4\over27\Theta^2}+{2\over81\Theta^3}\right),$$
is positive at small times, when $Q$ is sufficiently small, whereby
$p(\zeta)$ has three real roots $\zeta_i$. The polynomial has then
two local extrema at points of different signs, hence it has roots of both
signs. Since by Vi\`ete's theorem the product of the roots is negative,
two roots are positive and one is negative. Inequality \rf{poly} implies
\BE 2A+B\le\zeta_2(Q),\EE{bou}
where $\zeta_2$ is the intermediate root
(i.e., the smaller of the two positive roots).

We determine now the largest time $t_c$, for which bound \rf{bou} holds.
Differentiating the equation for the roots of polynomial \rf{poly} in $Q$, we find
$${\partial\zeta_i\over\partial Q}=-\left(\left.{\partial p\over
\partial\zeta}\right|_{\zeta=\zeta_i}\right)^{-1}.$$
Consequently, on increasing $Q$, the root $\zeta_2$ monotonically
increases and $\zeta_3$ monotonically decreases till the two roots collide (i.e.
$\zeta_2=\zeta_3$) at a critical value $Q=Q_c$. Then $t_c=Q_c/(2|\bom_0|_{0,\gamma}$).
For $Q=Q_c$ the discriminant $\Delta$ vanishes, which implies
\BE Q_c=\sqrt{\left({2\over3\Theta}+{32\over9}\right)^2+
{4\over27\Theta^2}+{2\over81\Theta^3}}-{2\over3\Theta}-{32\over9}.\EE{qc}
For small $\gamma$ and hence large $\Theta\sim1/\gamma$, we have
$Q_c=\Theta^{-2}/48+{\rm O}(\Theta^{-3})\sim\gamma^2$.

So far, inequality \rf{bou} for $t<t_c=Q_c/(2|\bom_0|_{0,\gamma})$ is just
an a priori bound. However, one can easily check that its derivation can be
repeated, without any changes, for partial sums of the series $A$ and $B$
having an arbitrary but equal number of terms. The bound \rf{bou} being uniform
in the number of terms left in the partial sums, convergence of the series
follows from standard arguments, as well as the validity of the bound \rf{bou}
for the infinite series $2A+B$. This remark concludes the proof.

Clearly, if the bound \rf{bou} holds true up to the time $t=t_c$,
then the series defining $A$ and $B$ are convergent in the complex time disk
$|t|<t_c$. Using the inequalities for the coefficients of the series for $2A+B$,
it is easy to show that, in this disk, the series \rf{xiser} is not just
a formal -- but a genuine -- solution to \rf{3DDisplacement1}--\rf{3DDisplacement2}.
By virtue of the proven \cite{BaFr} uniqueness of the solution in the said
H\"older class, any Lagrangian solution takes the form of the series \rf{xiser}.
This guarantees that, within this disk, the solution $\bxi$ and its gradient
$\nL\bxi$ are analytic functions in time.

Finally, we observe that the solution constructed to this point is based
on a Taylor expansion about $t=0$. However, we can restart from any time
$\widetilde{t}>0$, at which vorticity has the C$^{0,\widetilde{\gamma}}$
regularity for some $\widetilde{\gamma}>0$ (which can be smaller than
the H\"older exponent $\gamma$ for the initial data).
Thus, temporal analyticity persists on the interval
$[\widetilde{t},\widetilde{t}+\widetilde{t}_*)$ for a certain positive
$\widetilde{t}_*$. This process can be continued until the vorticity ceases
to have any C$^{0,\widetilde{\gamma}}$ regularity at some time $t_*$, which
may be finite or infinite. Hence trajectories of fluid particles remain real
analytic functions of the time as long as the vorticity preserves its H\"older
continuity. This completes our proof of analyticity in time of fluid particle
trajectories as long as a the velocity field is C$^{1,\gamma}$-continuous for
some $\gamma$ such that $0<\gamma<1$.

\section{Concluding remarks}\label{s:conclusion}

Examination of the structure of the proof of time analyticity
given in Section~\ref{s:hoelder} reveals that the essential ingredients are
(i) the fact that $|fg|_{0,\gamma} \le |f|_{0,\gamma} |g|_{0,\gamma}$
(algebra property) and (ii) the boundedness in the H\"older norm of
the operator $\nL_i\nL_j \nabla^{-2}$, where $\nabla^{-2}$ is the
inverse Lagrangian Laplacian acting on functions of zero spatial
mean. There are obviously many other function spaces for which similar
properties hold. This will be discussed elsewhere.

Of course, none of the results derived here from Cauchy's Lagrangian
formulation of the 3D incompressible Euler equation shed light on the
issue of whether or not there is finite-time blow up. It is
interesting to observe that, in Cauchy's formulation, all the
nonlinear terms vanish if the solutions depend only on one spatial
coordinate, because all the gradients appearing in the equations are
then parallel. Our bounds on the various nonlinear terms are rough as
they do not take into account the cancellation occurring in the vector
products, for instance, the possible strong depletion of nonlinearity
which can be associated to the development of almost one-dimensional
solutions. At present such issues can be investigated only
through very precise numerical simulations. By performing the
time-stepping in the Lagrangian domain one can then fully exploit
the analyticity result.

Here we have established results for the Lagrangian structure of
strong solutions, as long as they exist. Recently there has been significant
progress on weak solutions where the velocity has much less spatial
smoothness, e.g. H\"older continuity C$^{0,\gamma}$ with $\gamma$
close to the Kolmogorov--Onsager value $1/3$ \cite{camillolaszlo}. It
would be interesting to investigate the Lagrangian structure of such
weak solutions.

\section*{Acknowledgments}

We are grateful to C.~Bardos, Y.~Brenier, E.~Kuznetsov, T.~Matsumoto,
S.S.~Ray, L.~Sz\'ekelyhidi
and B.~Villone for extensive fruitful discussions. VZ was supported
by the grant 11-05-00167-a from the Russian foundation for basic research.
His several visits to the Observatoire de la C\^ote d'Azur (France) were
supported by the French Ministry of Higher Education and Research.


\begin{thebibliography}{99}
\bibitem{ohkitani}
K.~Ohkitani. koji-ohkitani.staff.shef.ac.uk/MAS411/files411/ch01.pdf
This is a chapter of a forthcoming book.
\end{thebibliography}
\end{document}